\documentclass[a4paper,10pt]{article}

\usepackage{setspace}
\usepackage[bottom]{footmisc}
\usepackage{color}
\usepackage{subfig}
\usepackage{hyperref}
\usepackage{cleveref}
\usepackage{multirow}
\usepackage{cite}
\usepackage{amsmath} 
\usepackage{mathtools}
\usepackage{graphicx}
\usepackage{epsfig}
\usepackage{epstopdf}
\usepackage{float}
\usepackage{titling}
\usepackage[margin=1.0in]{geometry}
\usepackage{etoolbox}
 
\title{A compact subcell WENO limiting strategy using immediate neighbors for Runge-Kutta Discontinuous Galerkin Methods}
\makeatletter
\renewcommand\@date{{%
  \vspace{-\baselineskip}%
  \large\centering
  \begin{tabular}{@{}c@{}}
    S R Siva Prasad Kochi\textsuperscript{1} \\
    \normalsize siva.ksr@gmail.com
  \end{tabular}%
  \quad and\quad
  \begin{tabular}{@{}c@{}}
    M Ramakrishna\textsuperscript{2} \\
    \normalsize krishna@ae.iitm.ac.in
  \end{tabular}

  \bigskip

  \textsuperscript{1}Doctoral Candidate, Dept. of Aerospace Engg., IIT Madras.\par
  \textsuperscript{2}Professor, Dept. of Aerospace Engg., IIT Madras.

  \bigskip

  \today
}}
\makeatother
\begin{document}

\maketitle

\begin{abstract}
 A compact subcell WENO (CSWENO) limiter is proposed for the solution of hyperbolic conservation laws with Discontinuous Galerkin Method which uses only the immediate neighbors of a given cell. These neighbors are divided into the required stencil for WENO reconstruction and an existing WENO limiting strategy is used. Accuracy tests and results for one-dimensional and two-dimensional Burgers equation and one-dimensional and two-dimensional Euler equations for Cartesian meshes are presented using this limiter. Comparisons with the parent WENO limiter are provided wherever appropriate and the performance of the current limiter is found to be slightly better than the parent WENO limiter for higher orders.
 
\end{abstract}

\section{Introduction}

In this paper, we look at the solution of hyperbolic conservation laws with the Runge Kutta Discontinuous Galerkin (RKDG) method \cite{cs1}. The main advantage of this method is that, in a given cell, to advance the solution in time, we need the information from only the immediate neighbors. For higher orders, to control spurious oscillations near discontinuities, a limiter is used and weighted essentially non oscillatory (WENO) limiters are preferred as they maintain the order of the scheme. However, standard WENO limiters need the information from neighbors of the neighboring cells and the advantage of RKDG method is lost.
\\
\\
\noindent Zhong and Shu \cite{zs} addressed this issue and used the whole DG polynomial in a given cell for WENO reconstruction using only the immediate neighbors. Dumbser et al \cite{dzls} used a different strategy where the target cell is divided into subcells and an \textit{a posteriori} limiting strategy is used based on the Multi-dimensional Optimal Order Detection (MOOD) approach. This gives very good results but it is quite expensive computationally. 
\\
\\
\noindent We propose a different strategy, where the immediate neighbors are divided into subcells based on the order of the scheme to get the required stencil, appropriate values are assigned to the new cells and the framework given in Qiu and Shu \cite{qs1} is used for limiting. We call the limiting strategy as compact subcell WENO limiter or CSWENO limiter in short.
\\
\\
\noindent The paper is organized as follows. We describe the formulation of the discontinuous Galerkin method used for all our results in \cref{sec:formulation}, the proposed limiting procedure is described in \cref{sec:newWENOLimiter} and finally the testing of the limiter and the results are described in \cref{sec:results} and we conclude the paper in \cref{sec:conc}.

\section{Formulation of Discontinuous Galerkin Method}\label{sec:formulation}

\noindent Consider the nonlinear scalar conservation law as given below

\begin{equation}\label{governEqn}
 \frac{\partial u}{\partial t} + \frac{\partial f(u)}{\partial x} = 0, \qquad \qquad x\in[x_{L},x_{R}]=\mathbf{D}
\end{equation}

\noindent with the initial condition

\begin{displaymath}
 u(x,0) = u_{0}(x),
\end{displaymath}

\noindent We look at solving \eqref{governEqn} using the Discontinuous Galerkin method. We approximate the domain $\mathbf{D}$ by $K$ non overlapping elements whose domain is given by $\mathbf{I}^{k}=[x_{l}^{k},x_{r}^{k}]$. We will approximate the local solution as a polynomial of order $N=N_{p}-1$, where $N_{p}$ is the number of degrees of freedom of the approximation. This is termed to be $\mathbf{P}^{N}$ based Discontinuous Galerkin method. The approximation is given as:

\begin{equation}\label{modalForm}
 u_{h}^{k}(x,t) = \sum_{n=0}^{N} \hat{u}^{k}_{n}(t)\psi_{n}^{k}(x) \qquad \forall x\in\mathbf{I}^{k}
\end{equation}

\noindent Here, $u_{h}^{k}(x,t)$ is the approximate local polynomial solution, $\psi_{n}^{k}(x)$ is the local polynomial basis of approximation and $\hat{u}^{k}_{n}(t)$ are the degrees of freedom.
\\
\noindent Similarly, we will also approximate the flux $f(u)$ in the domain $\mathbf{D}$ as given below:

\begin{equation}\label{fluxApprox}
 f_{h}^{k}(u_{h}^{k}) = \sum_{n=0}^{N} \hat{f}^{k}_{n}(t)\psi_{n}^{k}(x) \qquad \forall x\in\mathbf{I}^{k}
\end{equation}

\noindent We have used the orthonormalized Legendre polynomials as the local polynomial basis as suggested by Hesthaven et al \cite{hestha1}. The following affine mapping is employed.

\begin{equation}\label{affineMap}
 x(r) = x_{l}^{k} + \frac{1+r}{2}h^{k}, \qquad h^{k} = x_{r}^{k} - x_{l}^{k} \qquad \forall r\in\mathbf{I}=[-1,1]
\end{equation}

\noindent The degrees of freedom $\hat{u}^{k}_{n}$ can be advanced in time by the following scheme obtained from the weak form of the governing equation:

\begin{equation}\label{weakFormScheme}
 \frac{d}{dt}\hat{u}_{h}^{k} = (\mathbf{M}^{k})^{-1}(\mathbf{S}^{k})^{T} \hat{f}_{h}^{k}({u}_{h}^{k}) - (\mathbf{M}^{k})^{-1} (f^{*}|_{r_{N_{p}}}e_{N_{p}} - f^{*}|_{r_{1}}e_{1})
\end{equation}

\noindent Here, $\hat{u}_{h}^{k} = [\hat{u}_{0}^{k} \ldots \hat{u}_{N}^{k}]^{T}$, $e_{i}$ is a vector of dimension $N_{p}$ which has zero entries everywhere except at the $i$th location, and $\mathbf{M}^{k}$ is the local mass matrix which is given as:

\begin{equation}\label{massMatrix}
 \mathbf{M}^{k} = \left[M_{ij}^{k}\right] = \left[\int_{x_{l}^{k}}^{x_{r}^{k}} \psi_{i}^{k}(x) \psi_{j}^{k}(x) \text{dx}\right]
\end{equation}

\noindent and $\mathbf{S}^{k}$ is the local stiffness matrix which is given by:

\begin{equation}\label{stiffnessMatrix}
 \mathbf{S}^{k} = \left[S_{ij}^{k}\right] = \left[\int_{x_{l}^{k}}^{x_{r}^{k}} \psi_{i}^{k}(x) \frac{d\psi_{j}^{k}(x)}{dx} \text{dx}\right]
\end{equation}

\noindent Also, $f^{*}$ is the monotone numerical flux at the interface which is calculated using an exact or approximate Riemann solver. We have used the Lax-Friedrichs flux for all the test cases given below.

\section{Proposed limiting Procedure}\label{sec:newWENOLimiter}

\begin{figure}[htbp]
\begin{center}
\includegraphics[scale=1.0]{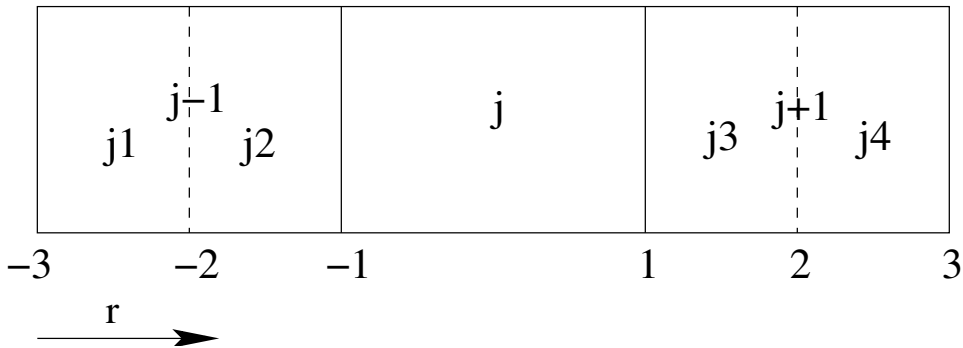}
\caption{New Stencil used for the WENO5 polynomial reconstruction for compact subcell WENO limiting ($r$ is given by the affine mapping defined in \eqref{affineMap})}
\label{fig1:newFivePointStencil}
\end{center}
\end{figure}

\noindent The common method for limiting in Discontinuous Galerkin method is:
\\
\noindent \textbf{1)} Identify the cells which need to be limited. They are known as troubled cells. \\
\noindent \textbf{2)} Replace the solution polynomial in the troubled cell with a new polynomial that is less oscillatory but with the same cell average and order of accuracy.
\\
\\
\noindent For the first step, we have used the KXRCF troubled cell indicator \cite{kxrcf} for all the calculations done in this paper as it is rated highly by Qiu and Shu in \cite{qs2} on the basis of it's performance in detecting the discontinuities in various test problems. The second step is where we do the limiting process. We will follow the method given in \cite{qs1}, where we use the finite volume WENO reconstruction procedure in the troubled cell, but with a new approach.
\\
\\
\noindent After identifying the troubled cells, we would like to reconstruct the moments, i.e., the values of $\hat{u}_{n}^{j}$ as given in \eqref{modalForm} for the troubled-cell $\mathbf{I}_{j}$ for $n=1,\ldots,N$. That is, we retain the cell average $\hat{u}_{0}^{j}$ and reconstruct all the other degrees of freedom. This is done as given in \cite{qs1} except for a modification which we suggest is that the stencil used for the reconstruction consists of only the immediate neighbors. For the $\mathbf{P}^{1}$ based DGM, we will need to use WENO3 reconstruction and we use the two-point Gauss quadrature rule to reconstruct the solution polynomial. Here, the procedure for limiting remains the same as in \cite{qs1} which is a standard WENO reconstruction using the cell averages. 
\\
\\
\noindent For the $\mathbf{P}^{2}$ based DGM, we have to use WENO5 reconstruction method. We use the four-point Gauss-Lobatto quadrature rule to reconstruct the solution polynomial, but instead of a stencil containing five cells as given in \cite{qs1}, we propose the following compact stencil. We take just three cells and divide the neighbors of the troubled-cell on both sides in half. This way, we get the new five point stencil required for the WENO5 reconstruction procedure as shown in Figure \ref{fig1:newFivePointStencil}. Here, the neighbors $\mathbf{I}_{j-1}$ and $\mathbf{I}_{j+1}$ are divided in half to obtain the new cells $\mathbf{I}_{j1}$, $\mathbf{I}_{j2}$, $\mathbf{I}_{j3}$ and $\mathbf{I}_{j4}$ respectively. Now, we use the DG polynomial for $\mathbf{I}_{j-1}$ and take its average in the new cells $\mathbf{I}_{j1}$ and $\mathbf{I}_{j2}$ to obtain the cell averages for the new cells. Similarly, we can calculate the cell average for the cells $\mathbf{I}_{j3}$ and $\mathbf{I}_{j4}$. We can use the same procedure to divide the neighboring cells for $\mathbf{P}^{N}$ based DGM for any $N$ based on the required number of quadrature points.

\noindent Note that the grid size is halved for the newly formed cells $\mathbf{I}_{j1}$, $\mathbf{I}_{j2}$, $\mathbf{I}_{j3}$ and $\mathbf{I}_{j4}$. Now, we will proceed as in \cite{qs1}. The step wise details of the WENO reconstruction are given below:
\\
\noindent \textbf{Step 1:} For the $\mathbf{P}^{N}$ based DGM, we need a Gauss or Gauss-Lobatto quadrature rule which is accurate to at least $O(h^{2N+2})$. This means that the order of accuracy for the WENO reconstruction has to be at least $2N+1$. Therefore, for the $\mathbf{P}^{1}$ based DGM, we need the two point Gauss quadrature points given by $r=-1/\sqrt{3}$ and $r=1/\sqrt{3}$. Similarly, for $\mathbf{P}^{2}$ based DGM, we need the four point Gauss-Lobatto quadrature points given by $r=-1$, $r=-1/\sqrt{5}$, $r=1/\sqrt{5}$ and $r=1$ and for $\mathbf{P}^{3}$ based DGM, we use the four point Gauss quadrature points given by $r=-\sqrt{525+70\sqrt{30}}/35$, $r=-\sqrt{525-70\sqrt{30}}/35$, $r=\sqrt{525-70\sqrt{30}}/35$ and $r=\sqrt{525+70\sqrt{30}}/35$.
\\
\noindent \textbf{Step 2:} The WENO reconstruction as given in \cite{shu1} is described now which is to be performed in the given troubled cell $\mathbf{I}_{j}$. We identify $N+1$ small stencils $S_{i}$, $i=0,\ldots,N$ such that $\mathbf{I}_{j}$ belongs to each of them. We set $S_{i} = \bigcup_{l=0}^{N}\mathbf{I}_{j+i-l}$ with the understanding that for the $\mathbf{P}^{2}$ based DGM, $\mathbf{I}_{j-2}$, $\mathbf{I}_{j-1}$, $\mathbf{I}_{j+1}$ and $\mathbf{I}_{j+2}$ are replaced by $\mathbf{I}_{j1}$, $\mathbf{I}_{j2}$, $\mathbf{I}_{j3}$ and $\mathbf{I}_{j4}$ respectively as explained above. We also have the larger stencil $\mathbf{T}=\bigcup_{i=0}^{N}S_{i}$ which contains all the cells from the smaller stencils $S_{i}$.
\\
\noindent Now, we have a polynomial of degree $N$, $p_{i}(x)$ corresponding to the stencil $S_{i}$ such that it's cell average in each of the cells of the stencil $S_{i}$ agrees with the given cell average of $u$. We also have a polynomial of degree $2N$ reconstruction denoted by $Q(x)$ associated with the larger stencil $\mathbf{T}$, such that the cell average of $Q(x)$ in each of the cells of the stencil $\mathbf{T}$ agrees with the cell average of $u$ for that cell. The details of the construction of $p_{i}(x)$ and $Q(x)$ are given in \cite{shu1}.
\\
\noindent \textbf{Step 3:} Next, we find the linear weights denoted by $\gamma_{0},\ldots,\gamma_{N}$, which satisfy

\begin{equation}\label{linearWeights}
 Q(x_{G}) = \sum_{i=0}^{N} \gamma_{i} p_{i}(x_{G})
\end{equation}

\noindent where $x_{G}$ is a Gauss or Gauss-Lobatto quadrature point. A set of linear weights for each of the quadrature points is obtained. The value of the functions $Q(x)$ and $p_{i}(x)$ for each $i$ can be written as a function of the cell average of each cell in the stencil. This is used in WENO reconstruction. 
For example, for the $\mathbf{P}^{2}$ based DGM, for the Gauss-Lobatto quadrature point $r=-1$, we have:

\begin{equation}
 \left[L^{2}_{i}\right] = \left[T^{3}_{ij}\right] \left[C^{2}_{j}\right]
\end{equation}

\noindent where

\begin{displaymath}
\left[L^{2}_{i}\right] =  \begin{bmatrix} p_{0}(x_{G}) & p_{1}(x_{G}) & p_{1}(x_{G}) & Q(x_{G}) \end{bmatrix}^{T}
\end{displaymath}

\begin{displaymath}
 \left[C^{2}_{j}\right] = \begin{bmatrix} \hat{u}_{0}^{j-2} & \hat{u}_{0}^{j-1} & \hat{u}_{0}^{j} & \hat{u}_{0}^{j+1} & \hat{u}_{0}^{j+2} \end{bmatrix}^{T}
\end{displaymath}

\noindent and

\begin{displaymath}
 \left[T^{3}_{ij}\right] = \begin{bmatrix} -\frac{1}{4} & \frac{13}{12} & \frac{1}{6} & 0 & 0 \rule{0pt}{2.5ex} \\ 0 & \rule{0pt}{2.5ex}  \frac{1}{2} & \rule{0pt}{2.5ex}  \frac{2}{3} & -\frac{1}{6} & 0 \\ 0 & 0 & \rule{0pt}{2.5ex} \frac{13}{6} & \rule{0pt}{2.5ex} -\frac{23}{12} & \rule{0pt}{2.5ex} -\frac{3}{4} \\ \rule{0pt}{2.5ex}-\frac{1}{10} & \rule{0pt}{2.5ex}\frac{21}{30} & \rule{0pt}{2.5ex}\frac{17}{30} & \rule{0pt}{2.5ex}-\frac{13}{60} & \rule{0pt}{2.5ex}\frac{1}{20} \end{bmatrix}
\end{displaymath}

\noindent The linear weights are given by

\begin{displaymath}
 \gamma_{0} = \frac{2}{5}, \qquad \qquad \gamma_{1} = \frac{24}{45}, \qquad \qquad \gamma_{2} = \frac{1}{15}
\end{displaymath}

\noindent For the Gauss-Lobatto quadrature point $r=-1/\sqrt{5}$, we have:

\begin{equation}
 \left[L^{2}_{i}\right] = \left[T^{4}_{ij}\right] \left[C^{2}_{j}\right]
\end{equation}

\noindent where

\begin{displaymath}
 \left[T^{4}_{ij}\right] = \begin{bmatrix} -\frac{3+6\sqrt{5}}{60} & \frac{5+18\sqrt{5}}{60} & \frac{58-12\sqrt{5}}{60} & 0 & 0 \\ 0 & \frac{2\sqrt{5}-1}{30} & \frac{16}{15} & -\frac{1+2\sqrt{5}}{30} & 0 \\ 0 & 0 & \frac{58+12\sqrt{5}}{60} & \frac{5-18\sqrt{5}}{60} & \frac{6\sqrt{5}-3}{60} \\ \frac{15-69\sqrt{5}}{3000} & \frac{63\sqrt{5}-29}{600} & \frac{163}{150} & -\frac{63\sqrt{5}+29}{600} & \frac{15+69\sqrt{5}}{3000} \end{bmatrix}
\end{displaymath}

\noindent The linear weights are given by

\begin{displaymath}
 \gamma_{0} = \frac{235-33\sqrt{5}}{950}, \qquad \qquad \gamma_{1} = \frac{48}{95}, \qquad \qquad \gamma_{2} = \frac{235+33\sqrt{5}}{950}
\end{displaymath}

\noindent Similarly, for the Gauss-Lobatto quadrature point $r=1$:

\begin{equation}
 \left[L^{2}_{i}\right] = \left[T^{5}_{ij}\right] \left[C^{2}_{j}\right]
\end{equation}

\noindent where

\begin{displaymath}
 \left[T^{5}_{ij}\right] = \begin{bmatrix} -\frac{3}{4} & -\frac{23}{12} & \frac{13}{6} & 0 & 0 \\ 0 & \rule{0pt}{2.5ex}-\frac{1}{6} & \rule{0pt}{2.5ex}\frac{2}{3} & \rule{0pt}{2.5ex}\frac{1}{2} & 0 \\ 0 & 0 & \rule{0pt}{2.5ex}\frac{1}{6} & \rule{0pt}{2.5ex}\frac{13}{12} & \rule{0pt}{2.5ex}-\frac{1}{4} \\ \rule{0pt}{2.5ex}\frac{1}{20} & \rule{0pt}{2.5ex}-\frac{13}{60} & \rule{0pt}{2.5ex}\frac{17}{30} & \rule{0pt}{2.5ex}\frac{21}{30} & \rule{0pt}{2.5ex}-\frac{1}{10} \end{bmatrix}
\end{displaymath}

\noindent The linear weights are given by

\begin{displaymath}
 \gamma_{0} = \frac{1}{15}, \qquad \qquad \gamma_{1} = \frac{24}{45}, \qquad \qquad \gamma_{2} = \frac{2}{5}
\end{displaymath}

\noindent Finally, for the Gauss-Lobatto quadrature point $r=1/\sqrt{5}$,

\begin{equation}
 \left[L^{2}_{i}\right] = \left[T^{6}_{ij}\right] \left[C^{2}_{j}\right]
\end{equation}

\noindent where

\begin{displaymath}
 \left[T^{6}_{ij}\right] = \begin{bmatrix} \frac{6\sqrt{5}-3}{60} & \frac{5-18\sqrt{5}}{60} & \frac{58+12\sqrt{5}}{60} & 0 & 0 \\ 0 & -\frac{1+2\sqrt{5}}{30} & \frac{16}{15} & \frac{2\sqrt{5}-1}{30} & 0 \\ 0 & 0 & \frac{58-12\sqrt{5}}{60} & \frac{5+18\sqrt{5}}{60} & -\frac{3+6\sqrt{5}}{60} \\ \frac{15+69\sqrt{5}}{3000} & -\frac{63\sqrt{5}+29}{600} & \frac{163}{150} & \frac{63\sqrt{5}-29}{600} & \frac{15-69\sqrt{5}}{3000} \end{bmatrix}
\end{displaymath}

\noindent The linear weights are given by

\begin{displaymath}
 \gamma_{0} = \frac{235+33\sqrt{5}}{950}, \qquad \qquad \gamma_{1} = \frac{48}{95}, \qquad \qquad \gamma_{2} = \frac{235-33\sqrt{5}}{950}
\end{displaymath}

\noindent \textbf{Step 4:} As given by \cite{shu1}, we compute the smoothness indicator for each stencil $S_{i}$:

\begin{equation}\label{smoothnessIndicator}
\beta_{i} = \sum_{l=1}^{N}\int_{I_{j}} \Delta x_{j}^{2l-1} \left(\frac{\partial^{l}}{\partial x^{l}}p_{i}(x)\right)^{2} dx
\end{equation}

%
%

\noindent For example, for the $\mathbf{P}^{2}$ based DGM, the smoothness indicators are given as:

\begin{displaymath}
 \beta_{0} = \frac{100(\hat{u}_{0}^{j})^{2}+\left(204\hat{u}_{0}^{j-2}-404\hat{u}_{0}^{j-1}\right)\hat{u}_{0}^{j}+433(\hat{u}_{0}^{j-1})^{2}-462\hat{u}_{0}^{j-1}\hat{u}_{0}^{j-2}+129(\hat{u}_{0}^{j-2})^{2}}{12}
\end{displaymath}

\begin{displaymath}
 \beta_{1} = \frac{43(\hat{u}_{0}^{j+1})^{2}+\left(70\hat{u}_{0}^{j-1}-156\hat{u}_{0}^{j}\right)\hat{u}_{0}^{j+1}+156(\hat{u}_{0}^{j})^{2}-156\hat{u}_{0}^{j-1}\hat{u}_{0}^{j}+43(\hat{u}_{0}^{j-1})^{2}}{9}
\end{displaymath}

\begin{displaymath}
 \beta_{2} = \frac{129(\hat{u}_{0}^{j+2})^{2}+\left(204\hat{u}_{0}^{j}-462\hat{u}_{0}^{j+1}\right)\hat{u}_{0}^{j+2}+433(\hat{u}_{0}^{j+1})^{2}-404\hat{u}_{0}^{j}\hat{u}_{0}^{j+1}+100(\hat{u}_{0}^{j})^{2}}{12}
\end{displaymath}

\noindent \textbf{Step 5:} We compute the nonlinear weights as given below:

\begin{equation}\label{nonLinearWeights}
\omega_{i}=\frac{\bar{\omega}_{i}}{\sum_{i}\bar{\omega}_{i}}, \quad \bar{\omega_{i}} = \frac{\gamma_{i}}{\sum_{i}(\epsilon + \beta_{i})^{2}}
\end{equation}

\noindent Here $\epsilon$ is a small number which is usually taken to be $10^{-6}$. The final WENO approximation is given by

\begin{equation}\label{wenoApprox}
u_{G} = \sum_{j=0}^{N}\omega_{i}p_{i}(x_{G})
\end{equation}

\noindent \textbf{Step 6:} We obtain the reconstructed degrees of freedom based on the reconstructed point values $u(x_{G})$ at the Gauss or Gauss-Lobatto quadrature points $x_{G}$ and a numerical integration as

\begin{equation}
\hat{u}_{i}^{j} = \Delta x_{j} \sum_{G} w_{G} u(x_{G}) \psi_{i}^{j}(x_{G}) \quad i=1,\ldots,N
\end{equation}

\noindent where $w_{G}$'s are the Gaussian quadrature weights for the points $x_{G}$. For a non-orthonormal basis, we define 
\begin{equation}
D_{i}^{j} = \Delta x_{j} \sum_{G} w_{G} u(x_{G}) \psi_{i}^{j}(x_{G}) \quad i=1,\ldots,N
\end{equation}

\noindent , $\mathbf{B}^{j} = \left[D_{1}^{j}-\hat{u}_{0}^{j}\left[M_{01}^{j}\right] \quad \ldots \quad  D_{N}^{j}-\hat{u}_{0}^{j}\left[M_{0N}^{j}\right]\right]^{T}$, $\mathbf{X}^{j} = \left[\hat{u}_{0}^{j} \ldots \hat{u}_{N}^{j}\right]^{T}$ and 

\begin{equation}
 \mathbf{A}^{j} = \left[A\right] =  \begin{bmatrix} \left[M_{11}^{j}\right] & \ldots & \left[M_{1N}^{j}\right] \\ \ldots & \ldots & \ldots \\ \rule{0pt}{2.5ex}  \left[M_{N1}^{j}\right] & \ldots & \rule{0pt}{2.5ex} \left[M_{NN}^{j}\right] \end{bmatrix}
\end{equation}

\noindent Here, the terms $\left[M_{mn}^{j}\right]$ are given by \eqref{massMatrix}. Then the reconstructed degrees of freedom are given by $\mathbf{X}^{j} = (\mathbf{A}^{j})^{-1}\mathbf{B}^{j}$. This will work for any polynomial basis. Now, we can get the reconstructed polynomial solution in $\mathbf{I}_{j}$ by \eqref{modalForm}. This completes the WENO limiting procedure. We call this limiting procedure the compact subcell WENO limiting or CSWENO limiting in short.
\\
\\
\noindent When you have a system of equations of the form \eqref{governEqn} to solve, in order to achieve better results, the limiter is used with a local characteristic field decomposition as explained in \cite{shu1}.
\\
\noindent For the two dimensional case, we reconstruct the values of the required function $u$ in the troubled cells at the tensor-product Gauss or Gauss-Lobatto quadrature points for the rectangular structured grids. For this, we get two polynomials, one each in $x$ and $y$ directions. Now, we apply the procedure given in \cite{zs} for time integration. For solving a system of equations, we use this with a local characteristic field decomposition.
\\
\\
\noindent Now, the semi-discrete scheme given in \eqref{weakFormScheme} along with the limiter is discretized in time by using the TVD Runge-Kutta time discretization introduced in \cite{shu}. We have used a third order TVD Runge-Kutta time discretization for all our calculations.

\section{Results}\label{sec:results}

\noindent In this section, we look at some of the results obtained to demonstrate the performance of the limiter (called the compact subcell WENO limiter or CSWENO limiter) described above. All the results are obtained using RKDG method and the CSWENO limiter with a third order TVD Runge-Kutta scheme for time integration unless otherwise specified.

\subsection{Accuracy Tests}\label{subsec:accuracyTest}

\noindent We test the accuracy of the schemes with the CSWENO limiter for scalar and system problems for both one-dimensional and two-dimensional test cases. We present the results of the accuracy tests using one and two-dimensional Burgers equations and one and two-dimensional nonlinear Euler equations. We used both uniform and non-uniform meshes for all the test cases. The non uniform meshes are obtained using a $10\%$ random perturbation of each node of the uniform mesh. We show only the results with nonuniform meshes as representative test cases.
\\
\\
\textbf{Example 1:} We solve the one dimensional nonlinear scalar inviscid Burgers equation:
\begin{equation}\label{1dBurgers}
 \frac{\partial u}{\partial t} + \frac{\partial (u^{2}/2)}{\partial x} = 0, \qquad \qquad x \in [0,2\pi]
\end{equation}
\noindent with the initial condition $u(x,0)=0.5+\sin x$, with periodic boundary conditions. The solution is smooth till $t=1.0$. The exact solution can be obtained using the Newton-Raphson method as given in \cite{heoc}. The errors and numerical orders of accuracy are calculated at $t=0.5$ and are presented in Table \ref{table:1}. We can see that the CSWENO limiter maintains the order and magnitude of accuracy of the original DG method.
\\
\\
\begin{table}
\centering
\resizebox{\textwidth}{!}{%
\begin{tabular}{|c|c|c|c|c|c|c|c|c|c|}
\hline
\multirow{2}{*}{} &  & \multicolumn{4}{|c|}{DG without limiter} & \multicolumn{4}{|c|}{DG with limiter} \\ \cline{2-10} 
 & $K$ & $L_{1}$ error & Order & $L_{\infty}$ error & Order & $L_{1}$ error & Order & $L_{\infty}$ error & Order \\ \hline
 \multirow{5}{*}{$\mathbf{P}^{1}$} & 20 & 5.32E-03 &  & 1.24E-02 &  & 9.16E-03 &  & 2.75E-02 &  \\ \cline{2-10}
  & 40 & 1.31E-03 & 2.02 & 3.22E-03 & 1.95 & 2.53E-03 & 1.92 & 6.08E-03 & 2.18 \\ \cline{2-10}
  & 80 & 3.18E-04 & 2.04 & 3.22E-03 & 1.95 & 6.52E-04 & 1.96 & 9.59E-04 & 2.66 \\ \cline{2-10}
  & 160 & 7.82E-05 & 2.02 & 2.13E-04 & 1.94 & 1.67E-04 & 1.96 & 2.62E-04 & 1.87 \\ \cline{2-10}
  & 320 & 1.91E-05 & 2.03 & 5.68E-05 & 1.91 & 4.45E-05 & 1.91 & 6.42E-05 & 2.03 \\ \hline
  \multirow{5}{*}{$\mathbf{P}^{2}$} & 20 & 3.12E-04 &  & 4.14E-03 &  & 4.89E-04 &  & 4.5E-03 & \\ \cline{2-10}
  & 40 & 4.61E-05 & 2.76 & 5.25E-04 & 2.98 & 6.12E-05 & 2.99 & 4.56E-04 & 2.94 \\ \cline{2-10}
  & 80 & 6.42E-06 & 2.84 & 6.89E-05 & 2.93 & 8.14E-06 & 2.91 & 6.32E-05 & 2.85 \\ \cline{2-10}
  & 160 & 9.24E-07 & 2.79 & 9.87E-06 & 2.8 & 1.15E-06 & 2.82 & 1.02E-05 & 2.63 \\ \cline{2-10}
  & 320 & 1.31E-07 & 2.82 & 1.68E-06 & 2.55 & 2.05E-07 & 2.49 & 2.47E-06 & 2.05 \\ \hline
  \multirow{5}{*}{$\mathbf{P}^{3}$} & 20 & 2.02E-05 &  & 4.13E-04 &  & 2.05E-05 &  & 4.16E-04 & \\ \cline{2-10}
  & 40 & 1.21E-06 & 4.06 & 3.42E-05 & 3.59 & 1.23E-06 & 4.06 & 3.43E-05 & 3.60 \\ \cline{2-10}
  & 80 & 7.62E-08 & 3.99 & 2.83E-06 & 3.60 & 7.67E-08 & 4.00 & 2.85E-06 & 3.59 \\ \cline{2-10}
  & 160 & 4.87E-09 & 3.97 & 2.81E-07 & 3.33 & 4.92E-09 & 3.96 & 2.81E-07 & 3.34 \\ \cline{2-10}
  & 320 & 2.98E-10 & 4.03 & 1.72E-08 & 4.03 & 3.06E-10 & 4.01 & 1.73E-08 & 4.02 \\ \hline
\end{tabular}}
\caption{1D Burgers equation with the initial condition $u(x,0)=0.5+\sin x$, with periodic boundary conditions, $t=0.5$, Nonuniform mesh with $K$ elements, $L_{1}$ and $L_{\infty}$ errors for $\mathbf{P}^{1}$, $\mathbf{P}^{2}$ and $\mathbf{P}^{3}$ based DG}
\label{table:1}
\end{table}

\noindent \textbf{Example 2:} We solve the two dimensional nonlinear scalar inviscid Burgers equation:
\begin{equation}\label{2dBurgers}
 \frac{\partial u}{\partial t} + \frac{\partial (u^{2}/2)}{\partial x} + \frac{\partial (u^{2}/2)}{\partial y} = 0, \qquad \qquad x,y \in [0,2\pi]
\end{equation}
\noindent with the initial condition $u(x,y,0)=0.5+\sin (x+y)$, with periodic boundary conditions. The solution is smooth till $t=0.5$. The exact solution is one-dimensional based on $x+y$ and can be calculated in a similar way to the one-dimensional problem. The errors and numerical orders of accuracy are calculated at $t=0.25$ and are presented in Table \ref{table:2}. Again, we can see that the CSWENO limiter maintains the order and magnitude of accuracy of the original DG method.
\\
\\
\begin{table}
\centering
\resizebox{\textwidth}{!}{%
\begin{tabular}{|c|c|c|c|c|c|c|c|c|c|}
\hline
\multirow{2}{*}{} &  & \multicolumn{4}{|c|}{DG without limiter} & \multicolumn{4}{|c|}{DG with limiter} \\ \cline{2-10} 
 & $K\times K$ & $L_{1}$ error & Order & $L_{\infty}$ error & Order & $L_{1}$ error & Order & $L_{\infty}$ error & Order \\ \hline
 \multirow{5}{*}{$\mathbf{P}^{1}$} & 20$\times$20 & 8.72E-03 &  & 1.56E-01 &  & 1.25E-02 &  & 2.9E-01 &  \\ \cline{2-10}
  & 40$\times$40 & 2.04E-03 & 2.10 & 4.93E-02 & 1.66 & 3.01E-03 & 2.05 & 7.41E-02 & 1.97 \\ \cline{2-10}
  & 80$\times$80 & 4.76E-04 & 2.10 & 1.37E-02 & 1.85 & 7.14E-04 & 2.08 & 1.96E-02 & 1.92 \\ \cline{2-10}
  & 160$\times$160 & 1.02E-04 & 2.22 & 3.89E-03 & 1.82 & 1.62E-04 & 2.15 & 5.34E-03 & 1.88 \\ \hline
  \multirow{5}{*}{$\mathbf{P}^{2}$} & 20$\times$20 & 9.1E-04 &  & 5.97E-02 &  & 9.2E-04 &  & 6.13E-02 &  \\ \cline{2-10}
  & 40$\times$40 & 1.13E-04 & 3.01 & 8.05E-03 & 2.89 & 1.24E-04 & 2.89 & 8.21E-03 & 2.9 \\ \cline{2-10}
  & 80$\times$80 & 1.76E-05 & 2.68 & 1.28E-03 & 2.65 & 1.96E-05 & 2.66 & 2.02E-03 & 2.02 \\ \cline{2-10}
  & 160$\times$160 & 2.5E-06 & 2.82 & 1.92E-04 & 2.74 & 2.68E-06 & 2.87 & 2.88E-04 & 2.81 \\ \hline
  \multirow{5}{*}{$\mathbf{P}^{3}$} & 20$\times$20 & 2.14E-04 &  & 9.71E-03 &  & 2.18E-04 &  & 9.92E-03 &  \\ \cline{2-10}
  & 40$\times$40 & 1.34E-05 & 4.00 & 6.17E-04 & 3.98 & 1.37E-05 & 3.99 & 6.31E-04 & 3.97 \\ \cline{2-10}
  & 80$\times$80 & 8.35E-07 & 4.00 & 4.05E-05 & 3.93 & 8.36E-07 & 4.03 & 4.11E-05 & 3.94 \\ \cline{2-10}
  & 160$\times$160 & 5.19E-08 & 4.01 & 2.13E-06 & 4.25 & 5.19E-08 & 4.01 & 2.13E-06 & 4.27 \\ \hline
\end{tabular}}
\caption{2D Burgers equation with the initial condition $u(x,0)=0.5+\sin (x+y)$, with periodic boundary conditions, $t=0.25$, Nonuniform mesh with $K\times K$ elements, $L_{1}$ and $L_{\infty}$ errors for $\mathbf{P}^{1}$, $\mathbf{P}^{2}$ and $\mathbf{P}^{3}$ based DG}
\label{table:2}
\end{table}

\noindent \textbf{Example 3:} We solve the two dimensional Euler equations:
\begin{equation}\label{2dEulerEquations}
\textbf{U}_{t} + \textbf{f(U)}_{x} + \textbf{g(U)}_{y} = 0
\end{equation}
\noindent where $\textbf{U} = (\rho, \rho u, \rho v, E)^{T}$, $\textbf{f(U)}=u\textbf{U} + (0, p, 0, pu)^{T}$ and $\textbf{g(U)}=v\textbf{U} + (0, 0, p, pv)^{T}$ with $p = (\gamma -1)(E-\frac{1}{2}\rho (u^{2}+v^{2}))$ and $\gamma = 1.4$. Here, $\rho$ is the density, $(u,v)$ is the velocity, $E$ is the total energy and $p$ is the pressure. The initial conditions are given by $\rho(x,y,0) = 1+0.2\sin (x+y)$, $u(x,y,0) = 0.7$, $v(x,y,0) = 0.3$ and $p(x,y,0) = 1.0$ and we use periodic boundary conditions. The exact solution is given by $\rho(x,y,0) = 1+0.2\sin (x+y-t)$, $u(x,y,0) = 0.7$, $v(x,y,0) = 0.3$ and $p(x,y,0) = 1.0$. The errors in density and numerical orders of accuracy are calculated at $t=2\pi$ and are presented in Table \ref{table:3}. Again, we can see that the CSWENO limiter maintains the order and magnitude of accuracy of the original DG method.
\\
\\
\begin{table}
\centering
\resizebox{\textwidth}{!}{%
\begin{tabular}{|c|c|c|c|c|c|c|c|c|c|}
\hline
\multirow{2}{*}{} &  & \multicolumn{4}{|c|}{DG without limiter} & \multicolumn{4}{|c|}{DG with limiter} \\ \cline{2-10} 
 & $K\times K$ & $L_{1}$ error & Order & $L_{\infty}$ error & Order & $L_{1}$ error & Order & $L_{\infty}$ error & Order \\ \hline
 \multirow{5}{*}{$\mathbf{P}^{1}$} & 20$\times$20 & 2.89E-03 &  & 7.42E-03 &  & 8.31E-03 &  & 3.24E-02 &  \\ \cline{2-10}
  & 40$\times$40 & 4.87E-04 & 2.57 & 2.11E-03 & 1.81 & 9.2E-04 & 3.18 & 7.52E-03 & 2.11 \\ \cline{2-10}
  & 80$\times$80 & 7.11E-05 & 2.78 & 6.89E-04 & 1.62 & 1.13E-04 & 3.03 & 9.13E-04 & 3.04 \\ \cline{2-10}
  & 160$\times$160 & 1.82E-05 & 1.97 & 2.16E-04 & 1.67 & 3.42E-05 & 1.72 & 2.45E-04 & 1.90 \\ \hline
  \multirow{5}{*}{$\mathbf{P}^{2}$} & 20$\times$20 & 9.88E-05 &  & 8.4E-04 &  & 1.25E-04 &  & 8.87E-04 &  \\ \cline{2-10}
  & 40$\times$40 & 1.25E-05 & 2.98 & 1.35E-04 & 2.64 & 1.64E-05 & 2.93 & 1.32E-04 & 2.75 \\ \cline{2-10}
  & 80$\times$80 & 1.54E-06 & 3.02 & 1.62E-05 & 3.06 & 2.02E-06 & 3.02 & 1.81E-05 & 2.87 \\ \cline{2-10}
  & 160$\times$160 & 1.84E-07 & 3.07 & 2.09E-06 & 2.95 & 2.31E-07 & 3.13 & 2.34E-06 & 2.95 \\ \hline
  \multirow{5}{*}{$\mathbf{P}^{3}$} & 20$\times$20 & 2.65E-06 &  & 6.87E-05 &  & 6.25E-06 &  & 7.01E-05 &  \\ \cline{2-10}
  & 40$\times$40 & 1.67E-07 & 3.99 & 4.42E-06 & 3.96 & 4.15E-07 & 3.91 & 6.24E-06 & 3.49 \\ \cline{2-10}
  & 80$\times$80 & 1.05E-08 & 3.99 & 2.85E-07 & 3.96 & 2.63E-08 & 3.98 & 4.32E-07 & 3.85 \\ \cline{2-10}
  & 160$\times$160 & 6.54E-10 & 4.01 & 1.65E-08 & 4.11 & 1.51E-09 & 4.12 & 3.12E-08 & 3.79 \\ \hline
\end{tabular}}
\caption{2D Euler equations with the initial condition $\rho(x,y,0) = 1+0.2\sin (x+y)$, $u(x,y,0) = 0.7$, $v(x,y,0) = 0.3$ and $p(x,y,0) = 1.0$, with periodic boundary conditions, $t=2\pi$, Nonuniform mesh with $K\times K$ elements, $L_{1}$ and $L_{\infty}$ errors for density with $\mathbf{P}^{1}$, $\mathbf{P}^{2}$ and $\mathbf{P}^{3}$ based DG}
\label{table:3}
\end{table}

\noindent \textbf{Example 4:} We again solve the two dimensional Euler equations given by \eqref{2dEulerEquations} for the Isentropic Euler Vortex problem suggested by Shu \cite{shu1} as a test case.The exact solution is given by: \\ $\rho = \left(1 -  \left(\frac{\gamma - 1}{16\gamma \pi^{2}}\right)\beta^{2} e^{2(1-r^{2})}\right)^{\frac{1}{\gamma-1}}$, $u = 1 - \beta e^{(1-r^{2})} \frac{y-y_{0}}{2\pi}$, $v = \beta e^{(1-r^{2})} \frac{x-x_{0}-t}{2\pi}$, and $p = \rho^{\gamma}$, where $r=\sqrt{(x-x_{0}-t)^{2}+(y-y_{0})^{2}}$, $x_{0}=5$, $y_{0}=0$, $\beta=5$ and $\gamma = 1.4$. We initialize with the exact solution at $t=0$ and use periodic boundary conditions at the edges of the domain. The errors in density and numerical orders of accuracy are calculated at $t=2$ and are presented in Table \ref{table:4}. Again, we can see that the CSWENO limiter maintains the order and magnitude of accuracy of the original DG method.
\\
\\
\begin{table}
\centering
\resizebox{\textwidth}{!}{%
\begin{tabular}{|c|c|c|c|c|c|c|c|c|c|}
\hline
\multirow{2}{*}{} &  & \multicolumn{4}{|c|}{DG without limiter} & \multicolumn{4}{|c|}{DG with limiter} \\ \cline{2-10} 
 & $K\times K$ & $L_{1}$ error & Order & $L_{\infty}$ error & Order & $L_{1}$ error & Order & $L_{\infty}$ error & Order \\ \hline
 \multirow{5}{*}{$\mathbf{P}^{1}$} & 20$\times$20 & 2.45E-03 &  & 1.89E-01 &  & 4.37E-03 &  & 2.81E-01 &  \\ \cline{2-10}
  & 40$\times$40 & 5.69E-04 & 2.11 & 4.71E-02 & 2.01 & 1.02E-03 & 2.10 & 7.18E-02 & 1.97 \\ \cline{2-10}
  & 80$\times$80 & 1.14E-04 & 2.32 & 1.29E-02 & 1.87 & 2.49E-04 & 2.03 & 1.85E-02 & 1.96 \\ \cline{2-10}
  & 160$\times$160 & 2.21E-05 & 2.37 & 3.34E-03 & 1.95 & 5.69E-05 & 2.13 & 4.23E-03 & 2.13 \\ \hline
  \multirow{5}{*}{$\mathbf{P}^{2}$} & 20$\times$20 & 6.43E-04 &  & 9.5E-02 &  & 2.82E-03 &  & 1.55E-01 &  \\ \cline{2-10}
  & 40$\times$40 & 7.87E-05 & 3.03 & 1.17E-02 & 3.02 & 2.98E-04 & 3.24 & 1.83E-02 & 3.08 \\ \cline{2-10}
  & 80$\times$80 & 8.92E-06 & 3.14 & 1.32E-03 & 3.15 & 3.76E-05 & 2.99 & 2.57E-03 & 2.83 \\ \cline{2-10}
  & 160$\times$160 & 1.17E-06 & 3.10 & 1.72E-04 & 2.94 & 4.18E-06 & 3.17 & 3.49E-04 & 2.88 \\ \hline
  \multirow{5}{*}{$\mathbf{P}^{3}$} & 20$\times$20 & 4.79E-06 &  & 8.35E-05 &  & 7.23E-06 &  & 8.37E-05 &  \\ \cline{2-10}
  & 40$\times$40 & 4.63E-07 & 3.96 & 5.98E-06 & 3.80 & 4.63E-07 & 3.97 & 6.33E-06 & 3.73 \\ \cline{2-10}
  & 80$\times$80 & 3.12E-08 & 3.99 & 4.16E-07 & 3.85 & 3.12E-08 & 3.89 & 4.24E-07 & 3.90 \\ \cline{2-10}
  & 160$\times$160 & 1.81E-09 & 4.04 & 2.72E-08 & 3.93 & 1.81E-09 & 4.11 & 3.25E-08 & 3.71 \\ \hline
\end{tabular}}
\caption{2D Euler equations for the Isentropic Vortex problem with periodic boundary conditions, $t=2$, Nonuniform mesh with $K\times K$ elements, $L_{1}$ and $L_{\infty}$ errors for density with $\mathbf{P}^{1}$, $\mathbf{P}^{2}$ and $\mathbf{P}^{3}$ based DG}
\label{table:4}
\end{table}

\subsection{Test Cases With Shocks}\label{subsec:TestCasesShocks}

\begin{figure}[htbp]
  \centering
  \subfloat[Solution of Burgers Equation with $\mathbf{P}^{2}$ based DGM]{\label{fig1:Burgers01}\includegraphics[width=0.45\textwidth]{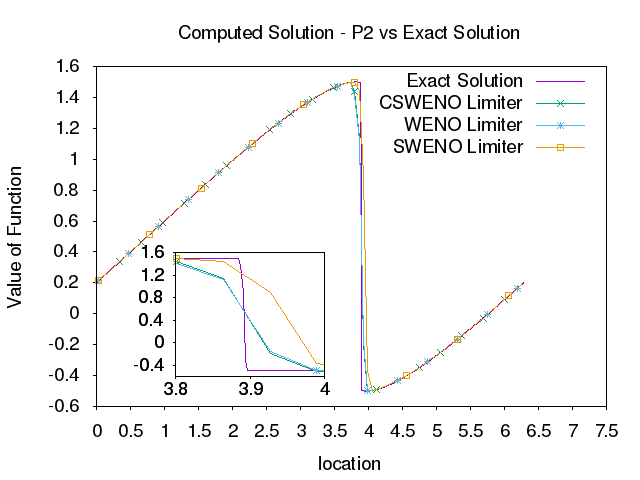}}
  \subfloat[Solution of Burgers Equation with $\mathbf{P}^{3}$ based DGM]{\label{fig2:Burgers02}\includegraphics[width=0.45\textwidth]{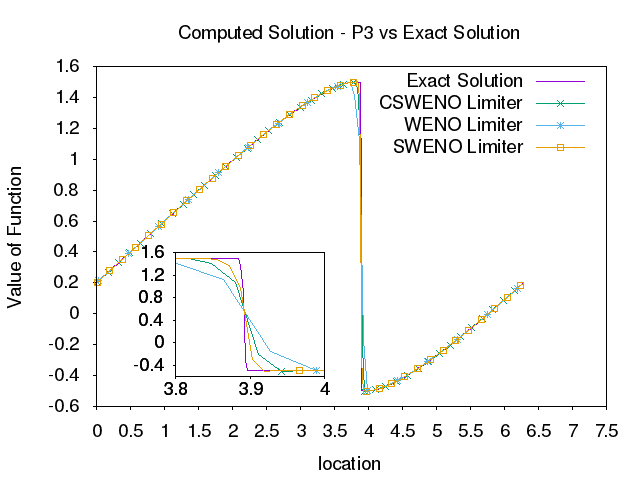}}\hfill
  \caption{Comparison of solutions of 1D Burgers Equation with $u(x,0)=0.5+\sin x$ at $t=1.5$ using 80 elements obtained with the CSWENO, WENO and SWENO limiters. Figures also include a zoomed in portion of the solution for better comparison}
  \label{fig:Burgers}
\end{figure}

\noindent We will now test the CSWENO (compact subcell WENO) limiter for problems with solutions having shocks. Here also, we have used both uniform and nonuniform meshes and obtained similar results. We will only show the results with uniform meshes. In all the results, we have compared the solution obtained with the CSWENO limiter against that obtained with the WENO limiter given in \cite{qs1} (also called the parent WENO limiter) and the simple WENO (also called SWENO) limiter described in \cite{zs}. As the described limiter is identical to WENO limiter for the $\mathbf{P}^{1}$ based DGM, we have only showed the comparisons for $\mathbf{P}^{2}$ and $\mathbf{P}^{3}$ based DGM. The WENO limiter uses a standard WENO reconstruction while the simple WENO limiter reconstructs the whole DG polynomial. In order to not clutter the comparison of the results, we have used fixed intervals between all the marked solution points and represented them only with lines otherwise.
\\
\\
\noindent \textbf{Example 5:} We solve the same nonlinear Burgers equation given in \eqref{1dBurgers} as in \textbf{Example 1} with the same initial condition $u(x,0)=0.5+\sin x$, with periodic boundary conditions. We now plot the results at $t=1.5$ when a shock has already appeared in the solution. The computed solution obtained at $t=1.5$ using 80 elements while using the CSWENO limiter for $\mathbf{P}^{2}$ and $\mathbf{P}^{3}$ based DGM is compared and plotted against the solution obtained using the WENO limiter, simple WENO limiter and the exact solution in Figures \ref{fig1:Burgers01} and \ref{fig2:Burgers02}. From the Figures, we can see that the performance of SWENO limiter is better than the CSWENO limiter and WENO limiter and the performance of CSWENO limiter is much better than the parent WENO limiter for both $\mathbf{P}^{2}$ based DGM and $\mathbf{P}^{3}$ based DGM.
\\
\\

\begin{figure}[htbp]
  \centering
  \subfloat[Solution of Buckley-Leverett Equation with $\mathbf{P}^{2}$ based DGM]{\label{fig1:BuckleyLeverett01}\includegraphics[width=0.45\textwidth]{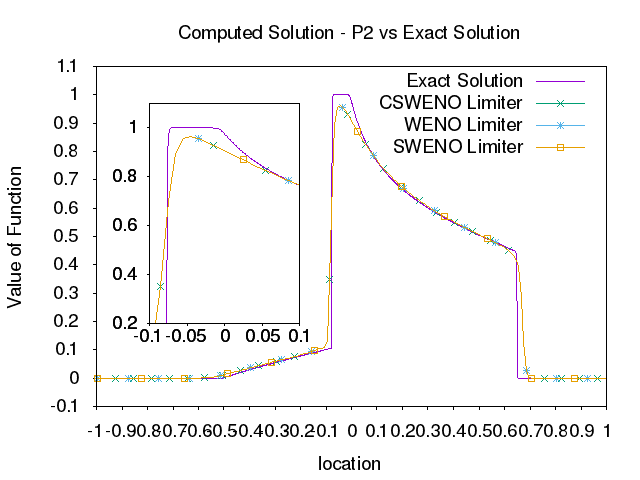}}
  \subfloat[Solution of Buckley-Leverett Equation with $\mathbf{P}^{3}$ based DGM]{\label{fig2:BuckleyLeverett02}\includegraphics[width=0.45\textwidth]{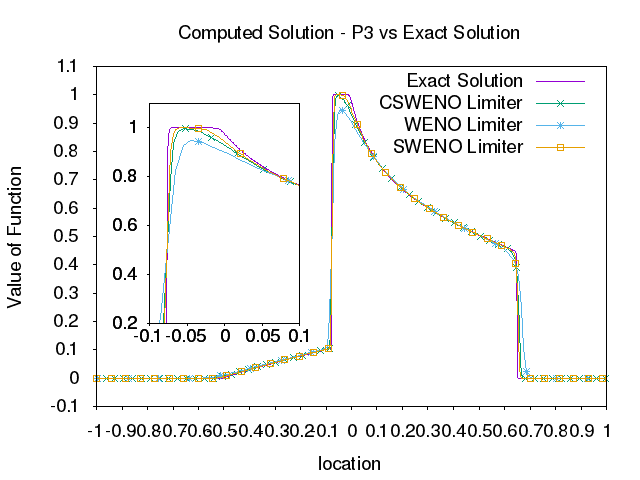}}\hfill
  \caption{Comparison of solutions of Buckley-Leverett Equation at $t=0.4$ using 80 elements obtained with the CSWENO, WENO and SWENO limiters. Figures also include a zoomed in portion of the solution for better comparison}
  \label{fig:BuckleyLeverett}
\end{figure}

\noindent \textbf{Example 6:} We solve the nonlinear nonconvex scalar Buckley-Leverett problem
\begin{equation}\label{1dBuckleyLeverett}
 \frac{\partial u}{\partial t} + \frac{\partial }{\partial x}\left(\frac{4u^{2}}{4u^{2}+(1-u)^{2}}\right) = 0, \qquad \qquad x \in [-1,1]
\end{equation}
with the initial condition $u=1$ for $-1/2\leq x \leq 0$ and $u=0$ everywhere else. The exact solution is a shock-rarefaction-contact discontinuity mixture. The computed solution obtained at $t=0.4$ using 80 elements while using the CSWENO limiter for $\mathbf{P}^{2}$ and $\mathbf{P}^{3}$ based DGM is compared and plotted against the solution obtained using the parent WENO limiter, simple WENO limiter and the exact solution in Figures \ref{fig1:BuckleyLeverett01} and \ref{fig2:BuckleyLeverett02}. From the Figures, we can see that the performance of all three limiters is quite similar for $\mathbf{P}^{2}$ based DGM. For $\mathbf{P}^{3}$ based DGM, performance of the SWENO limiter is slightly better than the CSWENO limiter and the parent WENO limiter and the performance of the CSWENO limiter is much better than the parent WENO limiter.
\\
\\

\begin{figure}[htbp]
  \centering
  \subfloat[Density solution of Sod Problem with $\mathbf{P}^{2}$ based DGM]{\label{fig1:Sod01}\includegraphics[width=0.45\textwidth]{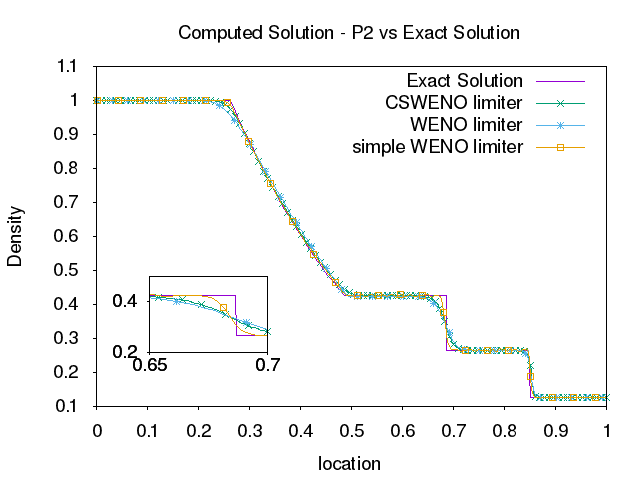}}
  \subfloat[Density solution of Sod Problem with $\mathbf{P}^{3}$ based DGM]{\label{fig2:Sod02}\includegraphics[width=0.45\textwidth]{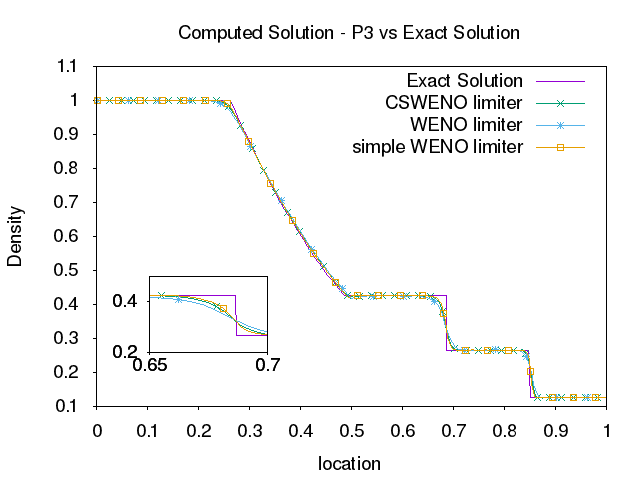}}\hfill
  \caption{Comparison of density solutions of Sod Problem at $t=0.2$ using 200 elements obtained with the CSWENO, WENO and SWENO limiters. Figures also include a zoomed in portion of the solution for better comparison}
  \label{fig:Sod}
\end{figure}

\noindent \textbf{Example 7:} We now consider a one-dimensional system particularly the Euler equations for an ideal gas given by

\begin{equation}\label{1dEulerEquations}
\textbf{U}_{t} + \textbf{f(U)}_{x} = 0
\end{equation}

\noindent where $\textbf{U}=(\rho, \rho u, E)^{T}$ and $\textbf{f(U)}=u\textbf{U} + (0, p, pu)^{T}$ with $p = (\gamma -1)(E-\frac{1}{2}\rho u^{2})$ and $\gamma = 1.4$. We will solve the Riemann problem in the domain $0\le x \le 1$ with the initial conditions given by Sod \cite{sod} as $(\rho_{L},u_{L},p_{L}) = (1,0,1)$ for $x<0.5$ and $(\rho_{R},u_{R},p_{R}) = (0.125,0,0.1)$ for $x\geq 0.5$. The computed solution for density obtained at $t=0.2$ using 200 grid points while using the CSWENO limiter for $\mathbf{P}^{2}$ and $\mathbf{P}^{3}$ based DGM is compared and plotted against the solution obtained using the parent WENO limiter, simple WENO limiter and the exact solution in Figures \ref{fig1:Sod01} and \ref{fig2:Sod02}. From the Figures, we can see that for $\mathbf{P}^{2}$ based DGM, performance of the SWENO limiter is slightly better than the CSWENO limiter and the parent WENO limiter and the performance of the CSWENO limiter is very similar to the parent WENO limiter. Also, the performance of all three limiters is quite similar for $\mathbf{P}^{3}$ based DGM.
\\
\\
%
%

\begin{figure}[htbp]
  \centering
  \subfloat[Density solution of Lax Problem with $\mathbf{P}^{2}$ based DGM]{\label{fig1:Lax01}\includegraphics[width=0.45\textwidth]{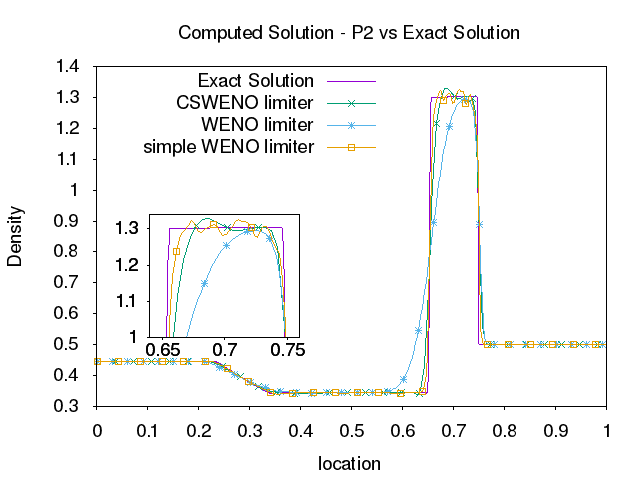}}
  \subfloat[Density solution of Lax Problem with $\mathbf{P}^{3}$ based DGM]{\label{fig2:Lax02}\includegraphics[width=0.45\textwidth]{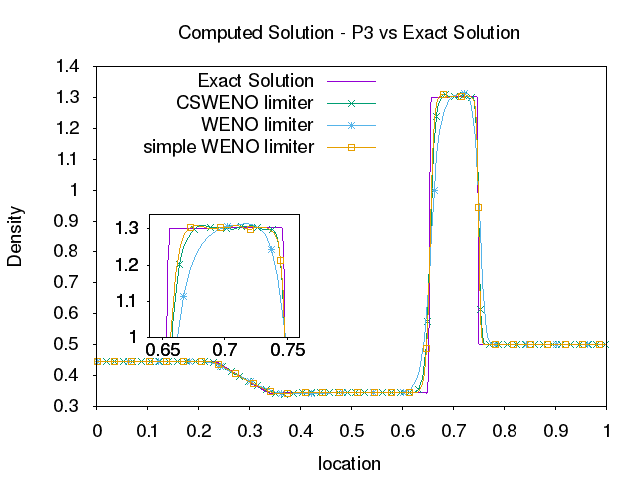}}\hfill
  \caption{Comparison of density solutions of Lax Problem at $t=0.1$ using 200 elements obtained with the CSWENO, WENO and SWENO limiters. Figures also include a zoomed in portion of the solution for better comparison}
  \label{fig:Lax}
\end{figure}

\noindent \textbf{Example 8:} We again solve the Euler equations for an ideal gas as given by \eqref{1dEulerEquations} for the Riemann problem in the domain $0\le x \le 1$ with the initial conditions given by Lax \cite{wc} as $(\rho_{L},u_{L},p_{L}) = (0.445,0.698,3.528)$ for $x<0.5$ and $(\rho_{R},u_{R},p_{R}) = (0.5,0,0.571)$ for $x\geq 0.5$. The computed solution for density obtained at $t=0.1$ using 200 elements while using the CSWENO limiter for $\mathbf{P}^{2}$ and $\mathbf{P}^{3}$ based DGM is compared and plotted against the solution obtained using the parent WENO limiter, simple WENO limiter and the exact solution in Figures \ref{fig1:Lax01} and \ref{fig2:Lax02}. Here, we can see that the simple WENO limiter performs slightly better than the parent WENO and the CSWENO limiter for both $\mathbf{P}^{2}$ based DGM and $\mathbf{P}^{3}$ based DGM and the CSWENO limiter performs quite better than the parent WENO limiter.
\\
\\

\begin{figure}[htbp]
  \centering
  \subfloat[Density solution of shock entropy wave problem with $\mathbf{P}^{2}$ based DGM]{\label{fig1:ShockDensityWave01}\includegraphics[width=0.45\textwidth]{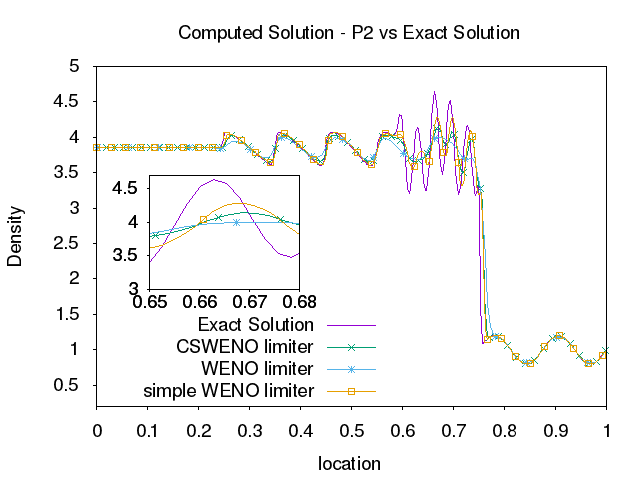}}
  \subfloat[Density solution of shock entropy wave problem with $\mathbf{P}^{3}$ based DGM]{\label{fig2:ShockDensityWave02}\includegraphics[width=0.45\textwidth]{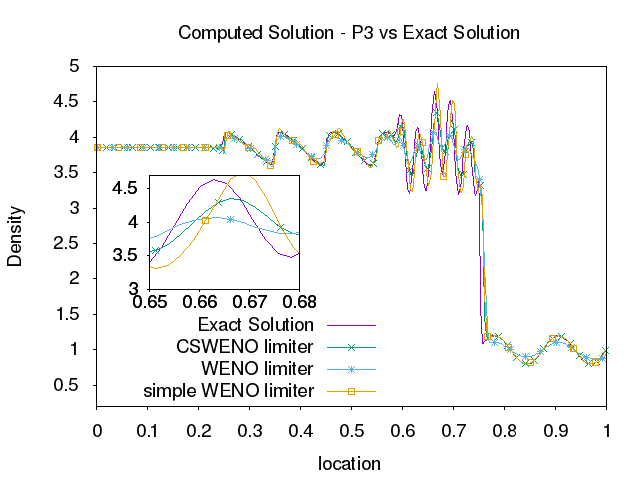}}\hfill
  \caption{Comparison of density solutions of Shock entropy wave Problem at $t=0.178$ using 200 elements obtained with the CSWENO, WENO and SWENO limiters. Figures also include a zoomed in portion of the solution for better comparison}
  \label{fig:ShockDensityWave}
\end{figure}

\noindent \textbf{Example 9:} We solve the problem of shock interaction with entropy waves as proposed in \cite{so2}. We solve the Euler equations \eqref{1dEulerEquations} with a moving shock interacting with sine waves in density in the domain $0\le x \le 1$ with the initial conditions given as $(\rho,u,p)=(3.857143,2.629369,10.333333)$ for $x < 0.125$ and $(\rho,u,p)=(1.0+0.2\sin(16\pi x),0,1)$ otherwise. The computed solution for density obtained at $t=0.178$s using 200 elements while using the CSWENO limiter for $\mathbf{P}^{2}$ and $\mathbf{P}^{3}$ based DGM is compared and plotted against the solution obtained using the parent WENO limiter, simple WENO limiter and the exact solution in Figures \ref{fig1:ShockDensityWave01} and  \ref{fig2:ShockDensityWave02}. From the Figures, we can see that for both $\mathbf{P}^{2}$ based DGM and $\mathbf{P}^{3}$ based DGM, the simple WENO limiter performs slightly better than the parent WENO and the CSWENO limiter. Also, the CSWENO limiter performs much better than the parent WENO limiter.
\\
\\

%

\begin{figure}[htbp]
  \centering
  \subfloat[Density solution of Blast Wave Problem with $\mathbf{P}^{2}$ based DGM]{\label{fig1:BlastWave01}\includegraphics[width=0.45\textwidth]{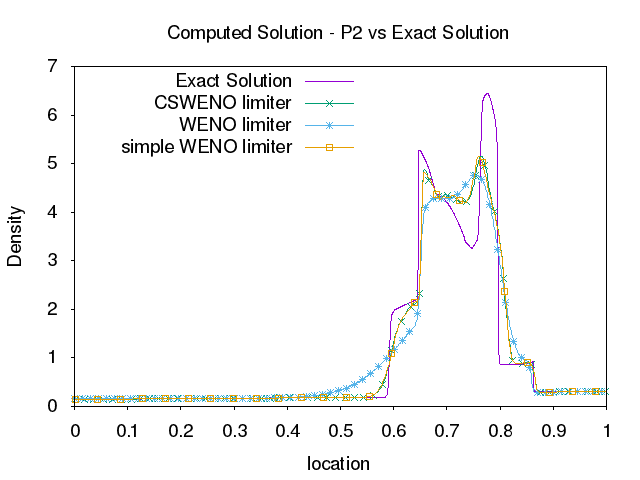}}
  \subfloat[Density solution of Blast Wave Problem with $\mathbf{P}^{3}$ based DGM]{\label{fig2:BlastWave02}\includegraphics[width=0.45\textwidth]{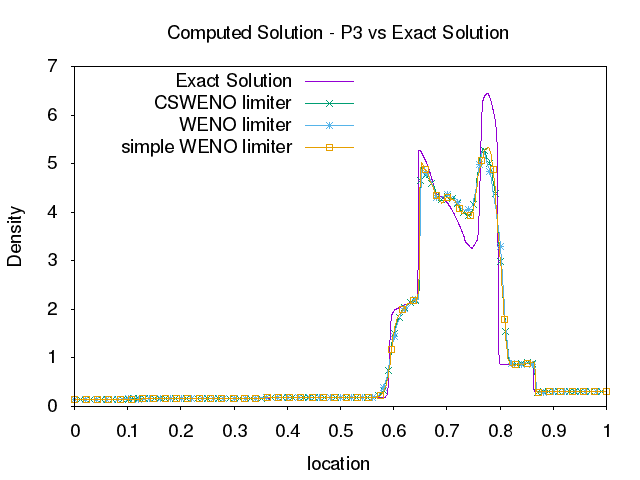}}\hfill
  \caption{Comparison of density solutions of Blast wave Problem at $t=0.038$ using 200 elements obtained with the CSWENO, WENO and SWENO limiters.}
  \label{fig:BlastWave}
\end{figure}

\noindent \textbf{Example 10:} We solve the interaction of blast waves of Euler equation as proposed in \cite{wc}. We solve the Euler equations \eqref{1dEulerEquations} in the domain $0\le x \le 1$ with the initial conditions given as $(\rho,u,p)=(1,0,1000)$ for $0 \leq x < 0.1$, $(\rho,u,p)=(1,0,0.01)$ for $0.1 \leq x < 0.9$ and $(\rho,u,p)=(1,0,100)$  otherwise with reflecting boundary conditions on both sides. The computed solution for density obtained at $t=0.038$s using 200 elements while using the CSWENO limiter for $\mathbf{P}^{2}$ and $\mathbf{P}^{3}$ based DGM is compared and plotted against the solution obtained using the parent WENO limiter, simple WENO limiter and the exact solution in Figures \ref{fig1:BlastWave01} and \ref{fig2:BlastWave02}. From the figures, for $\mathbf{P}^{2}$ based DGM, we can see that the performance of the CSWENO limiter and the simple WENO limiter is quite similar and both perform better than the WENO limiter. For $\mathbf{P}^{3}$ based DGM, the performance of all three limiters seems to be similar.
\\
\\

\begin{figure}[htbp]
  \centering
  \subfloat[Solution at t=0.2s with $\mathbf{P}^{2}$ based DGM]{\label{fig1:P1DGM}\includegraphics[width=0.9\textwidth]{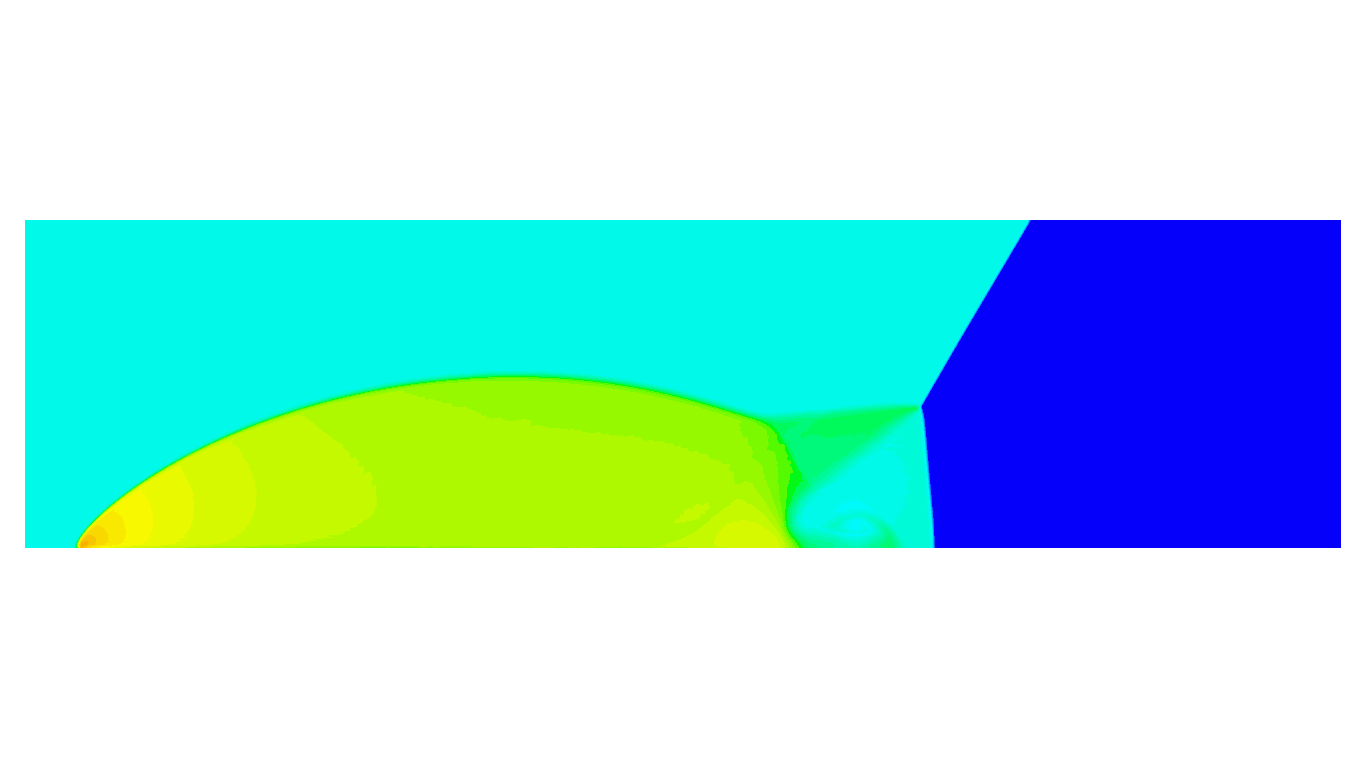}}\hfill
  \subfloat[Solution at t=0.2s with $\mathbf{P}^{3}$ based DGM]{\label{fig2:P2DGM}\includegraphics[width=0.9\textwidth]{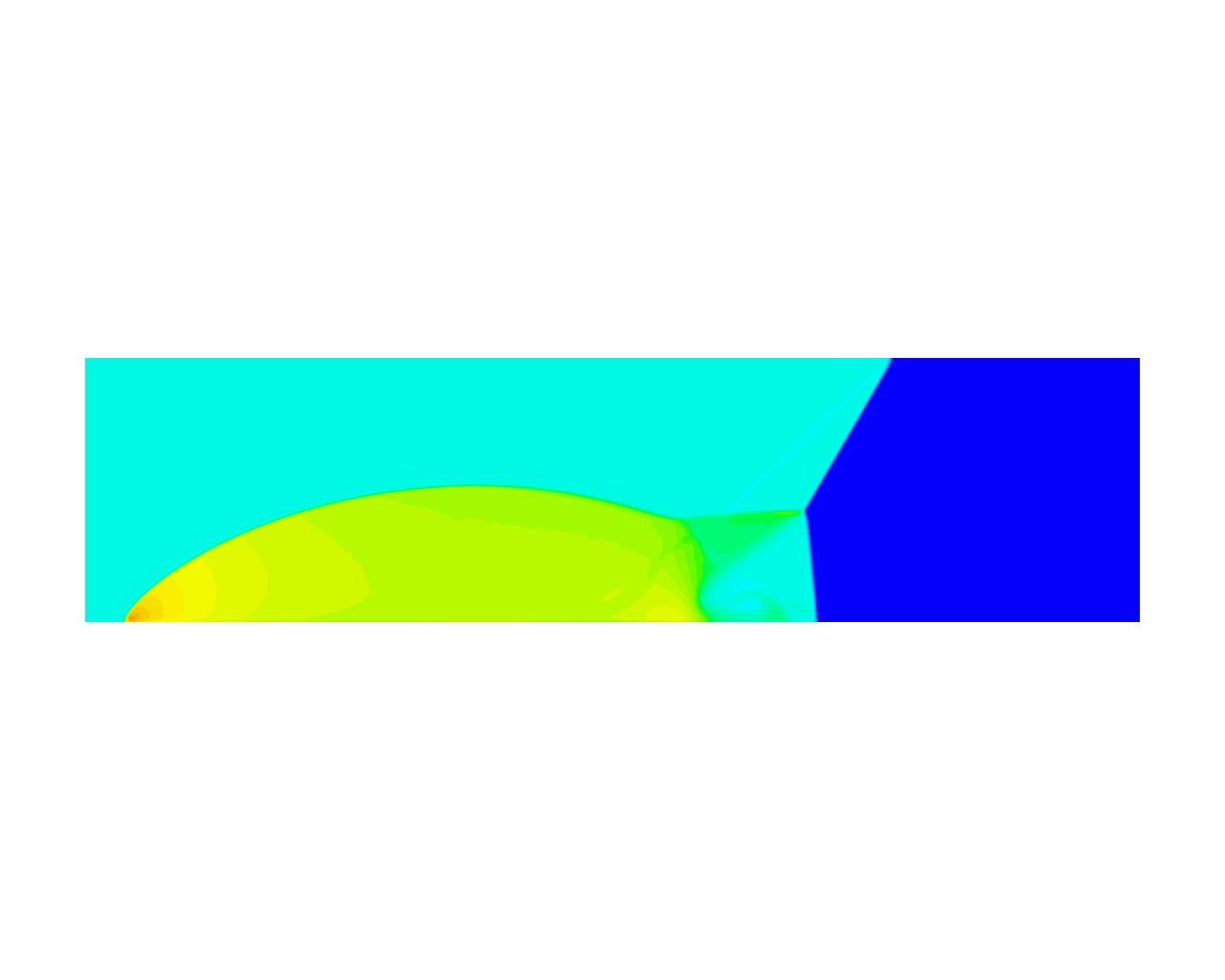}}\hfill
  \subfloat[Density Range]{\label{fig3:DensityRange}\includegraphics[width=0.9\textwidth]{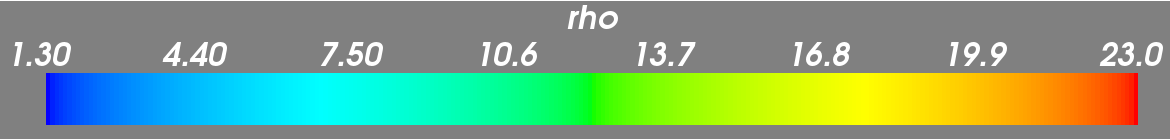}}
  \caption{Density variation for Double Mach reflection solution using the CSWENO limiter with $1920\times 480$ cells for $\mathbf{P}^{1}$(top) and $\mathbf{P}^{2}$(bottom) based DGM using 30 equally spaced contours}
  \label{fig8:newWENODMRSolution}
\end{figure}

\begin{figure}[htbp]
  \centering
  \subfloat[$960\times 240$ cells with $\mathbf{P}^{2}$ based DGM]{\label{fig3:P1DGM960new}\includegraphics[width=0.45\textwidth]{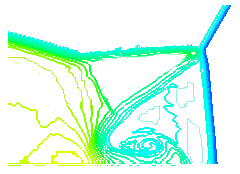}}
  \subfloat[$960\times 240$ cells with $\mathbf{P}^{3}$ based DGM]{\label{fig4:P2DGM960new}\includegraphics[width=0.45\textwidth]{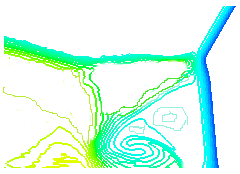}}\hfill
  \subfloat[$1920\times 480$ cells with $\mathbf{P}^{2}$ based DGM]{\label{fig5:P1DGM1920new}\includegraphics[width=0.45\textwidth]{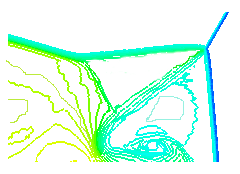}}
  \subfloat[$1920\times 480$ cells with $\mathbf{P}^{3}$ based DGM]{\label{fig6:P2DGM1920new}\includegraphics[width=0.45\textwidth]{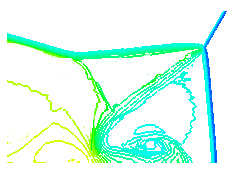}}\hfill
  \subfloat[Density Range]{\label{fig7:DensityRangenewWENO}\includegraphics[width=0.9\textwidth]{DensityRange.png}}
  \caption{Density variation for Double Mach reflection solution using the CSWENO limiter for two different mesh sizes in the region $[2,2.9]\times [0,0.6]$ using 30 equally spaced contours}
  \label{fig9:newWENODMRSolutionZoom}
\end{figure}

\begin{figure}[htbp]
  \centering
  \subfloat[$960\times 240$ cells with $\mathbf{P}^{2}$ based DGM]{\label{fig3:P1DGM960WENO}\includegraphics[width=0.45\textwidth]{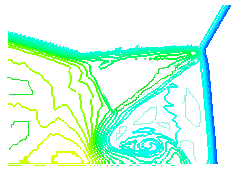}}
  \subfloat[$960\times 240$ cells with $\mathbf{P}^{3}$ based DGM]{\label{fig4:P2DGM960WENO}\includegraphics[width=0.45\textwidth]{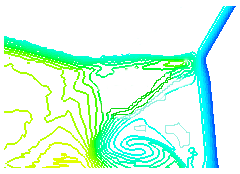}}\hfill
  \subfloat[$1920\times 480$ cells with $\mathbf{P}^{2}$ based DGM]{\label{fig5:P1DGM1920WENO}\includegraphics[width=0.45\textwidth]{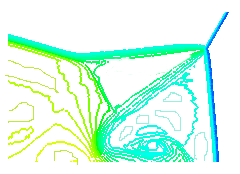}}
  \subfloat[$1920\times 480$ cells with $\mathbf{P}^{3}$ based DGM]{\label{fig6:P2DGM1920WENO}\includegraphics[width=0.45\textwidth]{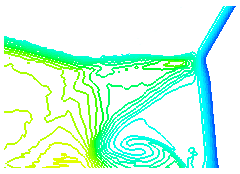}}\hfill
  \subfloat[Density Range]{\label{fig7:DensityRangeWENO}\includegraphics[width=0.9\textwidth]{DensityRange.png}}
  \caption{Density variation for Double Mach reflection solution using the WENO limiter for two different mesh sizes in the region $[2,2.9]\times [0,0.6]$ using 30 equally spaced contours}
  \label{fig10:WENODMRSolutionZoom}
\end{figure}

\begin{figure}[htbp]
  \centering
  \subfloat[$960\times 240$ cells with $\mathbf{P}^{2}$ based DGM]{\label{fig3:P1DGM960SWENO}\includegraphics[width=0.45\textwidth]{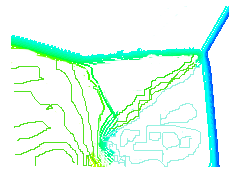}}
  \subfloat[$960\times 240$ cells with $\mathbf{P}^{3}$ based DGM]{\label{fig4:P2DGM960SWENO}\includegraphics[width=0.45\textwidth]{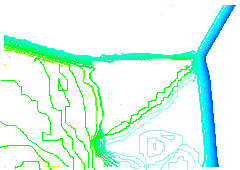}}\hfill
  \subfloat[$1920\times 480$ cells with $\mathbf{P}^{2}$ based DGM]{\label{fig5:P1DGM1920SWENO}\includegraphics[width=0.45\textwidth]{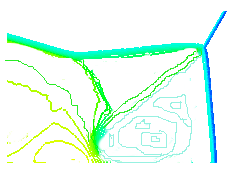}}
  \subfloat[$1920\times 480$ cells with $\mathbf{P}^{3}$ based DGM]{\label{fig6:P2DGM1920SWENO}\includegraphics[width=0.45\textwidth]{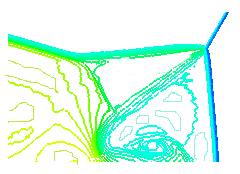}}\hfill
  \subfloat[Density Range]{\label{fig7:DensityRangeSWENO}\includegraphics[width=0.9\textwidth]{DensityRange.png}}
  \caption{Density variation for Double Mach reflection solution using the simple WENO limiter for two different mesh sizes in the region $[2,2.9]\times [0,0.6]$ using 30 equally spaced contours}
  \label{fig11:SWENODMRSolutionZoom}
\end{figure}

\noindent \textbf{Example 11:} As a test problem for the two-dimensional case, we solve the double Mach reflection problem which is given in \cite{wc}. We solve the two-dimensional Euler equations \eqref{2dEulerEquations} in the computational domain $[0,4]\times[0,1]$. Initially, right moving Mach 10 shock is positioned at $x=1/6,y=0$ and it makes an angle $60^{0}$ with the $x$-axis. For the bottom boundary, we impose the exact post shock conditions from $x=0$ to $x=1/6$ and for the rest of the $x$-axis, we use reflective boundary conditions. For the top boundary, we set conditions to describe the exact motion of a Mach 10 shock. We compute the solution upto time $t=0.2$ for two different uniform meshes with $960\times 240$ and $1920\times 480$ cells in each mesh. The full solution using the CSWENO limiter for the most refined mesh (containing $1920\times 480$ cells) for $\mathbf{P}^{2}$ and $\mathbf{P}^{3}$ based DGM has been shown in Figure \ref{fig8:newWENODMRSolution}. A zoom-in view of the density contours near the double Mach stem has been shown in Figures \ref{fig9:newWENODMRSolutionZoom}, \ref{fig10:WENODMRSolutionZoom} and \ref{fig11:SWENODMRSolutionZoom} using the CSWENO limiter, the parent WENO limiter and the SWENO limiter respectively for the two different mesh sizes. We have also added a positivity preserving component to the described limiter as given in \cite{zs3}.
\\
\\

We can see that the solution for all three limiters is quite well comparable to the solution obtained in \cite{wc}.

\section{Conclusions:}\label{sec:conc}

\noindent We have developed a different WENO limiting strategy for the solution of hyperbolic conservation laws using Discontinuous Galerkin method based on the limiter developed by Qiu and Shu \cite{qs1}. Here, we identify the troubled cells and use only the immediate neighbors by dividing them into subcells. This is different from the subcell limiting strategy of Dumbser et al \cite{dzls} which is much more accurate but quite complicated. These new cells are used for the reconstruction of the WENO polynomial. We termed this limiting procedure as the compact subcell WENO limiter (CSWENO limiter). We have tested the accuracy of this limiter using various standard test cases containing smooth solutions and calculating the numerical order of accuracy. We have also provided numerical results with shocks and compared the results obtained using this limiter with that obtained from the parent WENO limiter \cite{qs1} and the simple WENO limiter proposed by Zhong and Shu \cite{zs}. We can conclude from the results that the performance of the simple WENO limiter is better in most cases than the CSWENO limiter. Also, the CSWENO limiter performs better than the parent WENO limiter for most of the examples discussed. Also, the CSWENO limiter uses only the cell averages for the WENO reconstruction and uses a very compact stencil like the simple WENO limiter which is highly beneficial near the boundaries. Implementation of this limiter for unstructured meshes is ongoing.


\bibliographystyle{ieeetr}
\bibliography{references}

\begin{thebibliography}{10}

\bibitem{cs1}
B.~Cockburn and C.-W. Shu, ``The {Runge Kutta} local projection
  {$P^{1}$-discontinuous} {Galerkin} method for scalar conservation laws.,''
  {\em {$M^{2}AN$}}, vol.~25, pp.~337--361, 1991.

\bibitem{zs}
X.~Zhong and C.-W. Shu, ``A simple weighted essentially nonoscillatory limiter
  for {Runge-Kutta} discontinuous {Galerkin} methods.,'' {\em Journal of
  Computational Physics}, vol.~232, pp.~397--415, 2013.

\bibitem{dzls}
M.~Dumbser, O.~Zanotti, R.~Loubere, and S.~Diot, ``A posteriori subcell
  limiting of the discontinuous {Galerkin} finite element method for hyperbolic
  conservation laws.,'' {\em Journal of Computational Physics}, vol.~278,
  pp.~47--75, 2014.

\bibitem{qs1}
J.~Qiu and C.-W. Shu, ``{Runge-Kutta} discontinuous {Galerkin} method using
  {WENO} limiters.,'' {\em SIAM J. Sci. Comput.}, vol.~26, pp.~907--929, 2005.

\bibitem{hestha1}
J.~S. Hesthaven and T.~Warburton, {\em Nodal {Discontinuous} {Galerkin}
  {Methods:} {Algorithms,} {Analysis,} and {Applications}}.
\newblock Springer {New} {York}, 2008.

\bibitem{kxrcf}
L.~Krivodonova, J.~Xin, J.-F. Remacle, N.~Chevaugeon, and J.~Flaherty, ``Shock
  detection and limiting with discontinuous {Galerkin} methods for hyperbolic
  conservation laws.,'' {\em Appl. Numer. Math.}, vol.~48, pp.~323--338, 2004.

\bibitem{qs2}
J.~Qiu and C.-W. Shu, ``A comparison of troubled-cell indicators for
  {Runge-Kutta} discontinuous {Galerkin} methods using weighted essentially
  nonoscillatory limiters.,'' {\em SIAM J. Sci. Comput.}, vol.~27,
  pp.~995--1013, 2005.

\bibitem{shu1}
C.-W. Shu, ``Essentially non-oscillatory and weighted essentially
  non-oscillatory schemes for hyperbolic conservation laws.,'' {\em Lecture
  Notes in Mathematics, Springer}, vol.~1697, pp.~325--432, 1998.

\bibitem{shu}
C.-W. Shu, ``{TVD} time discretizations.,'' {\em SIAM J. Sci. Stat. Comput.},
  vol.~9, pp.~1073--1084, 1988.

\bibitem{heoc}
A.~Harten, B.~Engquist, S.~Osher, and S.~Chakravarthy, ``Uniformly high order
  accurate essentially {non-oscillatory} schemes, iii,'' {\em Journal of
  Computational Physics}, vol.~71, no.~2, pp.~231 -- 303, 1987.

\bibitem{sod}
G.~Sod, ``A survey of several finite difference methods for systems of
  nonlinear hyperbolic conservation laws.,'' {\em Journal of Computational
  Physics}, vol.~27, pp.~1--31, 1978.

\bibitem{wc}
P.~Woodward and P.~Colella, ``The numerical simulation of two-dimensional fluid
  flow with strong shocks.,'' {\em Journal of Computational Physics}, vol.~54,
  pp.~115--173, 1984.

\bibitem{so2}
C.-W. Shu and S.~Osher, ``Effective implementation of essentially
  non-oscillatory shock-capturing schemes, {II}.,'' {\em Journal of
  Computational Physics}, vol.~83, pp.~32--78, 1989.

\bibitem{zs3}
X.~Zhang and C.-W. Shu, ``On {positivity-preserving} high order discontinuous
  {Galerkin} schemes for compressible {Euler} equations on rectangular
  meshes.,'' {\em Journal of Computational Physics}, vol.~229, pp.~8918--8934,
  2010.

\end{thebibliography}

\end{document}